\documentclass{article}
\usepackage{latexsym}
\usepackage{leftidx}
\usepackage{amsmath,amsthm, color, pdfpages}
\usepackage{enumerate}
\usepackage{amssymb}
\usepackage{upgreek}
\usepackage[margin=1.1in,dvips]{geometry}

\def\f12{\frac 1 2}

\def\ga{\gamma}

\def\ep{\epsilon}

\def\La{\Lambda}

\def\Si{\Sigma}
\def\om{\omega}

\def\H{\mathcal{H}} 

\def\Lb{\underline{L}}

\def\Hb{\underline{\H}}

\def\pa{\partial}
\def\les{\lesssim}
\def\cL{{\mathcal L}}

\def\cD{\mathcal{D}} 
\def\f12{\frac 1 2}

\newcommand{\vol}{\textnormal{vol}}

\newcommand{\nabb}{\mbox{$\nabla \mkern-13mu /$\,}}

\newtheorem{thm}{Theorem}
\newtheorem{Thm}{Theorem}[section]

\newtheorem{Prop}{Proposition}[section]

\newtheorem{remark}{Remark}
\theoremstyle{definition}

\begin{document}

\title{Global behaviors of defocusing semilinear wave equations}
\date{}
\author{Shiwu Yang}
\maketitle
\begin{abstract}
In this paper, we investigate the global behaviors of solutions to defocusing semilinear wave equations in $\mathbb{R}^{1+d}$ with $d\geq 3$. We prove that in the energy space the solution verifies the integrated local energy decay estimates for the full range of energy subcritical and critical power. For the case when $p>1+\frac{2}{d-1}$, we derive a uniform weighted energy bound for the solution as well as inverse polynomial decay of the energy flux through hypersurfaces away from the light cone. As a consequence, the solution scatters in the energy space and in the critical Sobolev space for $p$ with an improved lower bound. This in particular extends the existing scattering results to higher dimensions without spherical symmetry.
\end{abstract}

\section{Introduction}
In this paper, we study the global asymptotic behaviors for solutions of  the following defocusing semilinear wave equation
\begin{equation}
  \label{eq:NLW:semi:alld}
  \Box\phi=|\phi|^{p-1}\phi,\quad \phi(0, x)=\phi_0(x),\quad \pa_t\phi(0, x)=\phi_1(x)
\end{equation}
with energy subcritical or energy critical power $1<p\leq \frac{d+2}{d-2}$ in $\mathbb{R}^{1+d}$.

This simple nonlinear model has draw extensive attention in the past decades. Existence of global classical $C^2$ solutions had been obtained early in \cite{Jorgens61:energysub:NLW:lowerd} with energy subcritical smooth nonlinearity in dimension $d=3$. Extensions and generalizations could be found for examples in \cite{Brenner79:globalregularity:NW}, \cite{Brenner81:globalregularity:d9}, \cite{Pecher76:NLW:global}, \cite{segal63:semigroup}, \cite{Vonwahl75:NW},
\cite{Strauss:NLW:decay}. These results aimed at showing the local boundedness of the solution with sufficiently regular initial data. However since the power $p$ is close to $1$ in higher dimensions, the above mentioned results only hold in lower dimensions $d\leq 9$. This calls for a global well-posedness result in the more natural energy space, which was addressed by Ginibre-Velo in \cite{velo85:global:sol:NLW}, \cite{Velo89:globalsolution:NLW} for the full energy subcritical case in all dimensions.
 This type of result is indeed a local existence result due to the conservation of energy and nothing too much could be said on the global and asymptotic behavior of the solutions.

 For the energy critical case, in dimension $d= 3$, Struwe in \cite{struwe88:NLW:d3:critical:symm} showed that the solution is globally in time if the data are spherically symmetric. Later Grillakis in \cite{Grillakis:NLW:cri:3d} removed this symmetry assumption and obtained the global regularity of the solution, that is, the solution is smooth if the initial data are. He also derived asymptotic pointwise decay estimates for the solution by using conformal transformation. This global regularity result was extended to higher dimensions up to $d\leq 9$ in
 \cite{Grillakis:NLW:cri:highD},  \cite{Shatah:criNLW:defocus:3D7}, \cite{Kapitanski94:NLW:n9:cri}.
On the other hand, in the energy space, Kapitanski \cite{Kapitanski90:NLW:globalweak} showed the existence and uniqueness of global weak solutions in the energy space and obtained partial regularity for the solution in \cite{Kapitanski90:NLW:H2regularity}. The global well-posedness in energy space was finally accomplished  by Shatah-Struwe in \cite{Shatah:criNLWdefocus:allD}.
A key observation leading to these global existence results is the non-concentration of energy. Bahouri-Shatah \cite{shatah98:timedecay:NLW:energycri} could even show that the potential part of energy decays to zero, which was applied to prove that the solution scatters to free linear wave by Bahouri-G\'{e}rard in \cite{Bahouri97:scattering:cri:NLW}.

Due to the lack of rigidity compared to the energy critical case, asymptotical properties for the solution in the energy subcritical case usually require additional restrictions on the initial data and the lower bound of the power $p$. In dimension $d=3$ with sufficiently smooth and localized initial data, pointwise time decay of the solution has been derived in
\cite{Strauss:NLW:decay}, \cite{vonWahl72:decay:NLW:super}, \cite{Roger:3DNLW:symmetry}, \cite{Bieli:3DNLW} for the superconformal case $p\geq \frac{d+3}{d-1}$. Weaker decay estimates were achieved in \cite{Pecher82:decay:3d}, \cite{Pecher82:NLW:2d} for part of subconformal case. These pointwise decay estimates in lower dimensions mainly relied on the approximate conservation of conformal energy (arising from the conformal symmetry of Minkowski space) as well as the representation formula of linear wave equations. However, such method fails in higher dimensions as the conserved energy is too weak to control the nonlinearity. Alternatively there have been plenty of literatures on the scattering theory of the solution, aiming at comparing the solution of nonlinear equation to that of linear equation at time infinity. A complete scattering theory consists of constructing a wave operator and proving the asymptotic completeness, that is, the solution behaves like linear solution at time infinity in a suitable function space (see detailed discussions in \cite{Velo87:decay:NLW}). Additional to the right function space for establishing such a theory is the lower bound of $p$ so that the nonlinearity decays sufficiently fast. Based on the approximate conservation of conformal energy,
 Ginibre-Velo \cite{Velo87:decay:NLW} obtained time decay for the solutions in the conformal energy space (weighted energy space with weights $(1+|x|^2)$) for the superconformal case, also see an alternative treatment in \cite{Baez:3DNLW:Groursat}. In lower dimensions $2\leq d\leq 4$, they also covered part of the subconformal case. These uniform time decay properties are sufficiently strong to establish the complete scattering theory, which was later extended to higher dimension $d\leq 6$ by Hidano \cite{Hidano03:scattering:NLW:56D}, \cite{Hidano:scattering:NLW}. However, it is not clear whether the lower bound given in these works is sharp or applies to even higher dimensions since the lower bound of $p$ for the asymptotic completeness is larger than that for existence of wave operator (see for example \cite{Hidano03:scattering:NLW:56D}).

One reason that $p$ is more restrictive for asymptotic completeness is that the norm of the chosen function space ($H^1$ with finite conformal energy) is too strong. The lower bound on $p$ could indeed be greatly improved for asymptotic completeness if the function space is enlarged to be the energy space $\dot{H}^1$, see
\cite{Pecher82:NLW:2d}, \cite{Pecher82:decay:SNW:4d}, \cite{Pecher82:decay:3d}. However these results still required that the initial data belong to the conformal energy space.

Another important intermediate function space to study the asymptotic completeness is the critical Sobolev space $\dot{H}^{s_p}$ motivated by the open problem that whether the nonlinear equation \eqref{eq:NLW:semi:alld} is globally well-posed and scattering in $\dot{H}^{s_p}$. Dodson \cite{Dodson:NLW:3d:p3}, \cite{Dodson:NLW:3d:p35} first gave an affirmative answer to this problem for the superconformal case $3\leq p<5$ in dimension $d=3$ under spherical symmetry. A conditional result, that is, the uniform boundedness of the critical Sobolev norm of the solution implies global existence and scattering, has been established in \cite{Dodson:NLW:3d:cond} without spherical symmetry. For data in some weighted energy space which belongs to the critical Sobolev space but contains the conformal energy space, Shen in \cite{Shen17:3DNLW:scatter} proved that the solution scatters in $\dot{H}^{s_p}$. However this result still requires spherical symmetry and is only for the superconformal case in dimension $d=3$.

The above mentioned results regarding the asymptotic decay properties of the solution do not hold for the one dimensional case $d=1$, $p>1$, as shown by Lindblad-Tao \cite{tao12:1d:NLW} that the solution exhibits a type of weak averaged decay estimate, which is clearly not shared by linear waves.
In particular, for $d=1$, the solution does not approach to linear ones as higher dimensional cases.

The aim of the present paper is to find new evidences that solutions to the energy subcritical and critical defocusing nonlinear wave equations behave like linear waves for $d\geq 3$ with initial data in some weighted energy space larger than the conformal energy space. There are several types of estimates that can characterize the global behavior of linear waves, of which the weakest version is the integrated local energy decay estimates. This type of estimate gives a uniform spacetime bound for the solution in terms of the initial energy and recently has been widely used to study linear waves on general Lorentzian manifolds, including black hole spacetimes, see for example \cite{dr3}, \cite{Tataru:localEdecay:Max}, \cite{Tataru13:localdecay}. We show that for the full range of energy subcritical and critical case $1<p\leq \frac{d+2}{d-2}$ solutions to \eqref{eq:NLW:semi:alld} verify the integrated local energy decay estimates, which in particular implies that the energy can not concentrated at a point. This fact is crucial to conclude the global well-posedness result for the energy critical equations.

Since linear waves travel along outgoing light cones, the energy flux through hypersurfaces away from the light cone decays in terms of the distance to the light cone as shown in \cite{newapp}. We demonstrate that for the case when $\frac{d+1}{d-1}<p\leq \frac{d+2}{d-2}$ quantitative energy flux decay estimates hold for solutions of \eqref{eq:NLW:semi:alld}. As a consequence, we have the uniform spacetime bound for the potential $|\phi|^{p+1}$, which leads to the scattering result in critical Sobolev space and in energy space with improved lower bound on $p$ than those in \cite{Velo87:decay:NLW}, \cite{Hidano03:scattering:NLW:56D}, \cite{Hidano:scattering:NLW} mentioned above. Moreover, our result applies to all higher dimensions $d\geq 3$ without spherical symmetry assumption, hence refining the scattering result of Shen in \cite{Shen17:3DNLW:scatter}.

\subsection{Statement of the main results}
We define some necessary notations.
The local coordinate system $(t, x)$, the associated polar local coordinate system $(t, r,
\om)$ with $r=|x|$, $\om=\frac{x}{|x|}$ as well as the null coordinates $u=\frac{t-r}{2}$, $v=\frac{t+r}{2}$ will be frequently used through out this paper.
We may use $u_+$, $v_+$ to denote $1+|u|$, $1+v$ respectively. For simplicity, we only consider the solution in the future $t\geq 0$. Same results hold for the past.
$\pa$ will be short for the full derivative $(\pa_t, \pa_{x^1}, \ldots, \pa_{x^d})$ and $\nabb$ means the covariant derivative on sphere with constant radius $r$ at fixed time $t$.

Let $\H_u$ be the outgoing null hypersurface $\{t-|x|=2u, |x|\geq 2\}$. Let $\Si_u=\H_u$ when $u<-1$ and
\begin{align*}
  \Si_{u}=\{|x|\leq 2, t=2u+2\}\cup\{\H_u\},\quad u\geq -1.
\end{align*}
This foliates the future of the Minkowski spacetime $\mathbb{R}^{1+d}$.
For $u\in \mathbb{R}$, let $\cD_{u}$ be the region consists of hypersurfaces $\H_{u'}$, $u'\leq u$ if $u< -1$. Otherwise $\cD_{u}$ stands for the future region inclosed by $\Si_{u}$. The energy flux through the hypersurface $\Si_u$ will be denoted as $E[\phi](\Si_u)$, that is
\begin{align*}
  E[\phi](\Si_u)=\int_{|x|\leq 2, t=2u+2}|\pa\phi|^2+\frac{2}{p+1}|\phi|^{p+1}dx+\int_{\H_u}(|\pa_v\phi|^2+|\nabb\phi|^2 +\frac{2}{p+1}|\phi|^{p+1}) r^{d-1}dvd\om
\end{align*}
when $u>-1$.

For $p>1$, denote
\[
s_p=\frac{d}{2}-\frac{2}{p-1}
\]
to be the critical exponent for the Sobolev norm $\dot{H}^{s_p}_x\times \dot{H}^{s_p-1}_x$ under scaling. The energy critical case corresponds to $s_p=1$ and the equation is conformally invariant when $s_p=\f12$.

Define the linear wave propagation operator $\mathbf{L}(t)$ as follows:
\[
\mathbf{L}(t)(\phi_0(x), \phi_1(x))=(\phi(t, x), \pa_t\phi(t, x)),
\]
in which $\phi$ solves the linear wave equation $\Box\phi=0$, $\phi(0, x)=\phi_0$, $\pa_t\phi(0, x)=\phi_1$.

We assume that the initial data are bounded in the following weighted energy space
\begin{align*}
  \mathcal{E}_{\ga_0}[\phi]=\int_{\mathbb{R}^d}(1+|x|)^{\ga_0}(|\nabla \phi_0|^2+| \phi_1|^2+\frac{2}{p+1}|\phi_0|^{p+1})dx
\end{align*}
for constant $\ga_0\in [0, 2]$. The case when $\ga_0=2$ corresponds to the conformal energy space and $\ga_0=0$ stands for the classical energy space. The situation considered in Shen's work \cite{Shen17:3DNLW:scatter} assumed that the initial data are spherically symmetric and bounded in $\mathcal{E}_{1+\ep}[\phi]$ for some $\ep>0$.

We are now ready to state our main results.
\begin{Thm}
  \label{thm:main}
  Consider the defocusing semilinear wave equation \eqref{eq:NLW:semi:alld} on $\mathbb{R}^{1+d}$, $d\geq 3$ with finite energy initial data $(\phi_0, \phi_1)$. Then the solution is globally in time and verifies the following asymptotical decay properties:
   \begin{itemize}
   \item For all $1<p\leq \frac{d+2}{d-2}$, the solution $\phi$ verifies an integrated local energy decay estimate
  \begin{equation}
  \label{eq:ILE:1p5}
    \iint_{\mathbb{R}^{1+d}} \frac{|\pa\phi|^2+|(1+r)^{-1}\phi|^2}{(1+r)^{1+\ep}}+\frac{|\phi|^{p+1}+|\nabb\phi|^2}{r} dxdt\leq C \mathcal{E}_{0}[\phi]
  \end{equation}
  for some constant $C$ depending only on $p$, $d$ and $\ep>0$.
  \item For all $\frac{d+1}{d-1}<p\leq \frac{d+2}{d-2}$ and $1<\ga_0<\min\{2, \f12 (p-1)(d-1)\}$, the solution $\phi$ satisfies the following quantitative inverse polynomial decay estimates
      \begin{align}
  \label{eq:ILE:cD:thm}
   E[\phi](\Si_{u})+\iint_{\cD_{u}} \frac{|\pa\phi|^2+|r_+^{-1}\phi|^2}{(1+r)^{1+\ep}}+\frac{|\phi|^{p+1}+|\nabb\phi|^2}{r} dxdt\leq C u_+^{-\ga_0}\mathcal{E}_{\ga_0}[\phi] ,\quad u\in \mathbb{R}
\end{align}
as well as the uniform $r$-weighted energy bound
\begin{align}
  \label{eq:PEW:thm}
  \iint_{\mathbb{R}^{1+d}} r^{\ga_0-d}|(\pa_t+\pa_r)(r^{\frac{d-1}{2}}\phi)|^2+r^{\ga_0-1}|\nabb\phi|^2+v_+^{\ga_0-\ep-1}|\phi|^{p+1}dxdt\leq C \mathcal{E}_{\ga_0}[\phi]
\end{align}
for some constant $C$ relying only on $p$, $d$, $\ep$ and $\ga_0$.
\item For all $p(d)=\frac{1+\sqrt{d^2+4d-4}}{d-1}<p< \frac{d+2}{d-2}$ and $$\max\{\frac{4}{p-1}-d+2, 1\}<\ga_0<\min\{\f12 (p-1)(d-1), 2\},$$
we have the uniform spacetime bound
\begin{equation}
  \label{eq:bd4:Lq:spacebd}
  \|\phi\|_{L^{\frac{(d+1)(p-1)}{2}}_{t, x}}\leq C  (p, d, \ga_0, \mathcal{E}_{\ga_0}[\phi])
\end{equation}
for some constant $C(p, d, \ga_0, \mathcal{E}_{\ga_0}[\phi]) $ depending on $p$, $d$, $\ga_0$ and the initial weighted energy norm $\mathcal{E}_{\ga_0}[\phi]$. As a consequence of the above uniform spacetime bound, the solution $\phi$ scatters to linear solutions in the Sobolev space $\dot H^s\times \dot H^{s-1}$ for all $s_p\leq s\leq 1$, that is, there exist pairs $\phi_0^{\pm}\in \dot H_x^{s_p}\cap \dot H_x^1$ and  $ \phi_1^{\pm}\in \dot H_x^{s_p-1}\cap L_x^2$ such that
\begin{align*}
      \lim\limits_{t\rightarrow\pm\infty}\|(\phi(t, x),\pa_t\phi(t,x))-\mathbf{L}(t)(\phi_0^{\pm}(x), \phi_1^{\pm}(x))\|_{\dot{H}_x^s\times \dot{H}_x^{s-1}}=0.
\end{align*}
  \end{itemize}
\end{Thm}
We give several remarks.
\begin{remark}
  As a consequence of the integrated local energy decay estimate, the potential energy can not concentrate at a point, that is, for any point $(t_0, x_0)$,
  \begin{align*}
    \liminf\limits_{t\rightarrow t_0} \int_{B_{t}}|\phi|^{p+1}dx =0,\quad B_t=\{(t, x)||x-x_0|\leq t_0-t\}.
  \end{align*}
  This estimate was crucial to conclude the global solution for \eqref{eq:NLW:semi:alld} with critical power $p$. To see that this estimate follows from the integrated local energy decay estimate \eqref{eq:ILE:1p5}, it suffices to move the origin to the point $(t_0, x_0)$. Then
  \begin{align*}
    \int_{t_1}^{t_0}\int_{B_{t}}(t_0-t)^{-1} |\phi|^{p+1} dxd t\leq \int_{t_1}^{t_0}\int_{B_{t}}\frac{|\phi|^{p+1}}{r}dxd t\leq C \mathcal{E}_0[\phi],
  \end{align*}
  which leads to the above claim.

  We point out here that using the same argument of Morawetz in \cite{mora1}, \cite{mora2}, one can also obtain the spacetime bound for the potential as well as the angular derivative of the solution (see \cite{Tao06:NDE:book}). The new ingredient of estimate \eqref{eq:ILE:1p5} is that it also controls the full derivative of the solution.
\end{remark}
\begin{remark}
  Both the integrated local energy decay estimate \eqref{eq:ILE:1p5} and the quantitative energy flux decay estimates \eqref{eq:ILE:cD:thm}, \eqref{eq:PEW:thm} hold for linear waves. In particular these estimates can be viewed as new evidences that solutions to \eqref{eq:NLW:semi:alld} asymptotically behave like linear solutions as long as $p>1+\frac{2}{d-1}$.
\end{remark}
\begin{remark}
  The lower bound $p^*(d)$ covered by Ginibre-Velo \cite{Velo87:decay:NLW} and Hidano \cite{Hidano03:scattering:NLW:56D} is the largest root of
  \[
(d-1)p^2-(d+2)p-1=0. 
\]
It can be checked that $p(d)<p^*(d)<\frac{d+3}{d-1}$. Moreover our scattering result holds for all dimension $d\geq 3$ while the above mentioned results only treated the lower dimension case $d\leq 6$.
\end{remark}

\begin{remark}
  A special case of the scattering result is when $3\leq p<5$, $d=3$, $\ga_0>1$, which has been investigated by Shen in \cite{Shen17:3DNLW:scatter} under spherical symmetry.
\end{remark}

\begin{remark}
  Estimates \eqref{eq:ILE:1p5}, \eqref{eq:ILE:cD:thm}, \eqref{eq:PEW:thm} also hold for $p>\frac{d+2}{d-2}$ if the solution $\phi$ is globally in time. However for supercritical equations, global well-posedness in energy space remains open. We restrict ourself to the energy critical and subcritical case is to guarantee the existence of global solution.
\end{remark}

\begin{remark}
  The lower bound $p(d)$ for the power $p$ arises for the existence of $\ga_0$ which records the decay rate of the initial data. The upper bound for $\ga_0$ is due to the $r$-weighted energy estimates while the lower bound for $\ga_0$ is used to control the nonlinearity for proving the scattering result. Although we have improved the lower bound for $p$, we conjecture that the sharp lower bound for the scattering result should be $p>1+\frac{4}{d}$. We note that
  \begin{align*}
    p(d)=\frac{1+\sqrt{d^2+4d-4}}{d-1}=1+\frac{4}{d}+\frac{8}{d^3}+O(\frac{1}{d^4}),
  \end{align*}
  which asymptotically approaches to $1+\frac{4}{d}$ to the third order in terms of the dimension $d$.
\end{remark}

\begin{remark}
  The uniform spacetime bound \eqref{eq:PEW:thm} will be used in the author's companion paper \cite{yang:NLW:ptdecay:3D} to show the improved pointwise decay estimates for the solution in dimension $d=3$.
\end{remark}

The proof for the integrated local energy estimates relies on the energy method by using the vector field $f(r)\pa_r$ as multiplier. The quantitative energy flux decay estimate as well as the $r$-weighted energy estimate follow by using the modified vector field method originally introduced by Dafermos-Rodnianski in \cite{newapp} for studying the linear waves on black holes spacetimes. The key of this new approach is a type of $r$-weighted energy estimates obtained by using the vector fields $r^{\ga}(\pa_t+\pa_r)$ as multipliers. Combined with the integrated local energy estimates as well as the energy conservation, a pigeonhole argument then leads to the inverse polynomial decay of the energy flux.
The main theorem can be viewed as a nonlinear application of this modified vector field method.

As a consequence of the $r$-weighted energy estimates, we also have the uniform weighted spacetime bound for the potential part $|\phi|^{p+1}$. By using a standard argument with the help of Strichartz estimates, we conclude the scattering result of the main theorem.

The plan of this paper is as follows: in Section 2, we will define some additional notations and review the vector field method. In Section 3, we obtain the integrated local energy decay estimates by using the vector field $f(r)\pa_r$ as multipliers. In Section 4, we apply the new approach to derive the $r$-weighted energy estimates as well as the quantitative energy flux decay. Then in the last section, we prove the uniform spacetime bound and conclude the scattering results.

\textbf{Acknowledgments.} The author would like to thank Mihalis Dafermos for enlightening and helpful discussions.  The work is partially supported by NSFC-11701017.

\section{Preliminaries and energy method}
\label{sec:notation}

Additional to the notations defined in the introduction, let $\{L, \Lb, e_1, e_2,
\ldots, e_{d-1}\}$ be  a null frame with
\[
L=\pa_v=\pa_t+\pa_r,\quad \Lb=\pa_u=\pa_t-\pa_r
\]
and $\{e_1, e_2, \ldots, e_{d-1}\}$ an orthonormal basis of the sphere with
constant radius $r$. At any fixed point $(t, x)$, we may choose this basis such that
\begin{equation}
 \label{eq:nullderiv}
 \begin{split}
&\nabla_{e_A}L=r^{-1}e_A,\quad \nabla_{e_A}\Lb=-r^{-1}e_A, \quad \nabla_{e_A}e_B=-\delta_{AB}r^{-1}\pa_r.
 \end{split}
\end{equation}
Here $\nabla$ is the covariant derivatives in Minkowski space and the index capital letter $A$ ranges from $1$ to $d-1$.

Similar to the outgoing null hypersurface $\H_u$, let $\Hb_v$ be the incoming null hypersurface $\{t+|x|=2v, |x|\geq 2\}$. We may also use the truncated ones $\H_u^{v_1, v_2}$, $\Hb_{v}^{u_1, u_2}$ defined as follows
\begin{align*}
  \H_u^{v_1, v_2}=\H_u\cap\{v_1\leq v\leq v_2\}, \quad
  \Hb_v^{u_1, u_2}=\Hb_v\cap\{u_1\leq u\leq u_2\}
\end{align*}
as well as the truncated foliation $\Si_u^v=\Si_u\cap\{t+|x|\leq 2v\}$.
The exterior region will be referred as $\{(t, x)| u=\frac{t-|x|}{2}\leq -1, t\geq 0\}$  while the interior region is $\{(t, x)| u\geq -1, t\geq 0\}$.

We now can define the domains that would be used in this paper. Let $\mathcal{D}_{u_1, u_2}^{v}$  be the domain
\[
\mathcal{D}_{u_1, u_2}^{v}:=\{\cup_{u_1\leq u\leq u_2}\Si_{u}^v\}
\]
for $u_1\leq u_2$ and $u_2\leq -1$ or $u_1\geq -1$. In particular these regions are located in the exterior region or in the interior region. We may omit the index $v$ when $v=\infty$.

For simplicity, for integrals in this paper, we will omit the volume form unless it is specified. More precisely we will use
\begin{align*}
  \int_{\mathcal{D}}f,\quad \int_{\H} f, \quad \int_{\Hb}f, \quad \int_{\{t=constant\}} f
\end{align*}
to be short for
\begin{align*}
  \int_{\mathcal{D}}f dxdt, \quad \int_{\H} f 2r^{d-1}dvd\om , \quad \int_{\Hb}f 2r^{d-1}dud\om, \quad \int_{\{t=constant\}} f dx
\end{align*}
respectively. Here $\om$ are the standard coordinates of unit sphere.

Through out this paper, we make a convention that $A\les B$ means that there exists a constant $C$, depending possibly on $p$, $d$, $\ga_0$ and some small positive constant $\ep$ such that $A\leq CB$.
\bigskip

Now we review the energy method for wave equations.
Define the associated energy momentum tensor for the scalar field $\phi$
\begin{align*}
  T[\phi]_{\mu\nu}=\pa_{\mu}\phi\pa_{\nu}\phi-\f12 m_{\mu\nu}(\pa^\ga \phi \pa_\ga\phi+\frac{2}{p+1} |\phi|^{p+1}),
\end{align*}
where $m_{\mu\nu}$ is the flat Minkowski metric on $\mathbb{R}^{1+d}$. Then we can compute that
\begin{align*}
  \pa^\mu T[\phi]_{\mu\nu}=&(\Box\phi-  |\phi|^{p-1}\phi)\pa_\nu\phi.
\end{align*}
Now for any vector fields $X$, $Y$ and any function $\chi$, define the current
\begin{equation*}
J^{X, Y, \chi}_\mu[\phi]=T[\phi]_{\mu\nu}X^\nu -
\f12\pa_{\mu}\chi \cdot|\phi|^2 + \f12 \chi\pa_{\mu}|\phi|^2+Y_\mu.
\end{equation*}
Then for solution $\phi$ of equation \eqref{eq:NLW:semi:alld}, we have the energy identity
\begin{equation}
\label{eq:energy:id}
\iint_{\mathcal{D}}\pa^\mu  J^{X,Y,\chi}_\mu[\phi] d\vol =\iint_{\mathcal{D}}div(Y)+ T[\phi]^{\mu\nu}\pi^X_{\mu\nu}+
\chi \pa_\mu\phi\pa^\mu\phi -\f12\Box\chi\cdot|\phi|^2 +\chi \phi\Box\phi d\vol
\end{equation}
for any domain $\mathcal{D}$ in $\mathbb{R}^{1+d}$. Here $\pi^X=\f12 \cL_X m$  is the deformation tensor of the metric $m$ along the vector field $X$.

\section{The integrated local energy decay estimates}

In this section, we derive the integrated local energy decay estimates for solution $\phi$ of \eqref{eq:NLW:semi:alld}. Take
\begin{align*}
  X=f(r)\pa_r,\quad \chi=\frac{d-1}{2}r^{-1}f,\quad Y=0
  \end{align*}
for some function $f(r)$ of the radius $r=|x|$. Under the coordinate system $(t, x)$, we compute the non-vanishing components of the deformation tensor
\begin{align*}
  \pi^X_{ij}=(f'-r^{-1}f)\om_i\om_j+r^{-1}fm_{ij}.
\end{align*}
We therefore can compute that
\begin{align*}
  & T[\phi]^{\mu\nu}\pi^X_{\mu\nu}+
\chi \pa_\mu\phi \pa^\mu\phi  +\chi\phi\Box\phi-\f12 \Box\chi |\phi|^2\\
&=(f'-r^{-1}f)(|\pa_r\phi|^2-\f12 (\pa^\ga\phi\pa_\ga\phi+\frac{2}{p+1}|\phi|^{p+1}))+r^{-1}f (|\nabla\phi|^2-\frac{d}{2}(\pa^\ga\phi\pa_\ga\phi+\frac{2}{p+1}|\phi|^{p+1}))\\
&\quad +
\chi \pa_\mu\phi \pa^\mu\phi  +\chi|\phi|^{p+1}-\f12 \Box\chi |\phi|^2\\
&=\f12 f'(|\pa_t\phi|^2+|\pa_r\phi|^2)+(r^{-1}f-\f12 f')|\nabb\phi|^2+\frac{(p-1)(d-1)r^{-1}f-2f'}{2(p+1)}|\phi|^{p+1}-\f12 \Box\chi |\phi|^2.
\end{align*}
Now take
\[
f=2\delta_p^{-1}+1-r_+^{-\ep},\quad r_+=1+r,\quad \delta_p=(p-1)(d-1)
\]
for some small positive constant $\ep<\f12$. We then can compute that
\begin{align*}
  &f'=\ep r_+^{-1-\ep},\quad r^{-1}f-f'=\frac{r_+^{1+\ep}-1-(1+\ep)r}{r_+^{1+\ep}}+2\delta_p^{-1}r^{-1}\geq 2\delta_p^{-1}r^{-1},\\
  &-\Box \chi=\frac{d-1}{2}r^{-1}\left(\ep(1+\ep)r_+^{-2-\ep}+(d-3)r^{-1}(r^{-1}f-f')\right).
\end{align*}
This shows that
\begin{align*}
  (p-1)(d-1)r^{-1}f-2f'\geq \delta_p (f'+2\delta_p^{-1}r^{-1})-2f'=2(r^{-1}-\ep r_+^{-1-\ep})\geq r^{-1}.
\end{align*}
Thus the above calculations show that
\begin{align*}
 \frac{|\pa\phi|^2+|r_+^{-1}\phi|^2}{(1+r)^{1+\ep}}+\frac{|\phi|^{p+1}+|\nabb\phi|^2}{r}\les  T[\phi]^{\mu\nu}\pi^X_{\mu\nu}+
\chi \pa_\mu\phi \pa^\mu\phi  +\chi\phi\Box\phi-\f12 \Box\chi |\phi|^2.
\end{align*}
Here the implicit constant relies only on $p$, $d$ and $\ep$.

Apply the energy identity \eqref{eq:energy:id} to the domain $\cD_{u_1, u_2}^{v_0}$ for $-1\leq u_1\leq u_2$, $v_0\geq 2+u_2$ or $ u_1\leq u_2<-1$, $v_0=-u_1$. The above computation leads to
\begin{align}
\label{eq:ILE:0:cD}
  \iint_{\cD_{u_1, u_2}^{v_0}} \frac{|\pa\phi|^2+|r_+^{-1}\phi|^2}{(1+r)^{1+\ep}}+\frac{|\phi|^{p+1}+|\nabb\phi|^2}{r}\les |\int_{\pa \cD_{u_1, u_2}^{v_0} } i_{J^{X, Y, \chi}[\phi]}d\vol|.
\end{align}
Similarly for the domain bounded by the initial hypersurface and the $t$-constant slice, we obtain that
\begin{align}
\label{eq:ILE:0:Rd}
  \int_0^t\int_{\mathbb{R}^d}\frac{|\pa\phi|^2+|r_+^{-1}\phi|^2}{(1+r)^{1+\ep}}+\frac{|\phi|^{p+1}+|\nabb\phi|^2}{r} \les \left| \int_{\mathbb{R}^d} J^{X, Y, \chi}[\phi]^0dx|_{0}^t \right|.
\end{align}
As $\cD_{u_1, u_2}^{v_0}$ is bounded by the hypersurfaces $\Si_{u_1}^v$, $\Si_{u_2}^v$ and $\Hb_v^{u_1, u_2}$, we estimate the boundary terms on the outgoing null hypersurface $\H_u$, the incoming null hypersurface $\Hb_v$ and the $t$-constant hypersurface. For the outgoing null hypersurface $\H_u$, we can show that
\begin{align*}
|i_{J^{X, Y, \chi}[\phi]}d\vol |&= | T[\phi]_{L \nu}X^\nu -
\f12 L\chi |\phi|^2 + \f12 \chi\cdot L |\phi|^2 |  r^{d-1}dv d\om\\
&=\f12|  f (|L\phi|^2  -|\nabb\phi|^2-\frac{2}{p+1}|\phi|^{p+1})-
 L\chi |\phi|^2 +   2\chi \phi\cdot L \phi |  r^{d-1}dv d\om\\
 &\les (  f (|L\phi|^2  +|\nabb\phi|^2+\frac{2}{p+1}|\phi|^{p+1})+(\chi^2+|\chi'|)|\phi|^2
 )  r^{d-1}dv d\om.
\end{align*}
By definition, we can bound that
\begin{align*}
  |f|\leq 2\delta_p^{-1}+1,\quad \chi^2\leq r_+^{-2}((d-1)\delta_p^{-1}+\frac{d-1}{2})^2,\quad |\chi'|\leq (2\delta_p^{-1}+1)r_+^{-2}.
\end{align*}
Therefore on the outgoing null hypersurface $\H_u$, we can bound that
\begin{align*}
|i_{J^{X, Y, \chi}[\phi]}d\vol | \les (  |L\phi|^2  +|\nabb\phi|^2+\frac{2}{p+1}|\phi|^{p+1}+r_+^{-2}|\phi|^2
 )  r^{d-1}dv d\om.
\end{align*}
The first three terms in the above integrand are exactly the energy flux for solution $\phi$ of the semilinear wave equation \eqref{eq:NLW:semi:alld} and the last term could be bounded by using Hardy's inequality.

Next on the incoming null hypersurface $\Hb_v$, similarly we first can estimate that
\begin{align*}
  |i_{J^{X, Y, \chi}[\phi]}d\vol|&= | T[\phi]_{\Lb \nu}X^\nu -
\f12 \Lb\chi |\phi|^2 + \f12 \chi\cdot \Lb |\phi|^2  |  r^{d-1}du d\om\\
&=\f12| f (|\nabb\phi|^2+\frac{2}{p+1}|\phi|^{p+1}-|\Lb\phi|^2) +
 \chi' |\phi|^2 +   2\chi \phi\cdot \Lb \phi  |  r^{d-1}du d\om\\
 &\les (|\Lb\phi|^2+|\nabb\phi|^2+\frac{2}{p+1}|\phi|^{p+1}+r_+^{-2}|\phi|^2)r^{d-1}du d\om.
\end{align*}
Here we used the above bounds for $f$, $\chi$ and $\chi'$.

On the constant $t$-hypersurface, we estimate that
\begin{align*}
  |i_{J^{X, Y,\chi}[\phi]}d\vol|&=| T[\phi]_{0 \nu}X^\nu -
\f12 \pa_t\chi |\phi|^2 + \f12 \chi\cdot \pa_t |\phi|^2 |  dx\\
&= | f\pa_t\phi\pa_r\phi  +   \chi  \phi \pa_t \phi|  dx\\
&\les (|\pa\phi|^2+\frac{2}{p+1}|\phi|^{p+1}+r_+^{-2}|\phi|^2)dx.
\end{align*}
Now for hypersurface $\Si$ in Minkowski space, denote
\begin{align*}
  E[\phi](\Si)=2\int_{\Si}i_{J^{\pa_t, 0, 0}[\phi]}d\vol.
\end{align*}
We then can compute that
\begin{align*}
  E[\phi](\H_u)&=\int_{\H_u}(|L\phi|^2+|\nabb\phi|^2+\frac{2}{p+1}|\phi|^{p+1})r^{d-1}dvd\om,\\
  E[\phi](\Hb_v)&=\int_{\Hb_v}(|\Lb\phi|^2+|\nabb\phi|^2+\frac{2}{p+1}|\phi|^{p+1})r^{d-1}dud\om,\\
  E[\phi](\mathbb{R}^d)&=\int_{\mathbb{R}^d}|\pa\phi|^2+\frac{2}{p+1}|\phi|^{p+1}dx.
\end{align*}
Thus applying the energy identity \eqref{eq:energy:id} to the vector fields $X=\pa_t$, $Y=0$ and function $\chi=0$, we obtain the energy conservation
\begin{align*}
  E[\phi](\Si_{u_1}^v)&=E[\phi](\Si_{u_2}^v)+E[\phi](\Hb_v^{u_1, u_2}),\\
  E[\phi](\{t\}\times \mathbb{R}^d)&=E[\phi](\{0\}\times \mathbb{R}^d)
\end{align*}
for all $t\geq 0$ and $-1\leq u_1<u_2$ as well as the energy conservation in the exterior region
\begin{align*}
  E[\phi](\H_{u_1}^{-u_2})+E[\phi](\Hb_{-u_2}^{u_1, u_2})=E[\phi](\{t=0, |u_1|\leq \frac{|x|}{2}\leq |u_2|\})
\end{align*}
for all $u_2\leq u_1\leq -1$. Letting $u_2\rightarrow -\infty$ and by definition of $\mathcal{E}_{\ga}[\phi]$, we obtain the energy flux decay in the exterior region
\begin{equation}
  \label{eq:Eflux:decay:ex}
  E[\phi](\H_u)\leq C u_+^{-\ga}\mathcal{E}_{\ga}[\phi],\quad \forall u\leq -1,\quad \ga\geq 0.
\end{equation}
This in particular proves the energy flux decay of the solution in the exterior region.

Now in view of the integrated local energy estimates \eqref{eq:ILE:0:cD} and by using Hardy's inequality (see for example \cite{yang1}, \cite{dr3}), we can show that
\begin{align*}
&\iint_{\cD_{u_1, u_2}^{-u_1}} \frac{|\pa\phi|^2+|r_+^{-1}\phi|^2}{(1+r)^{1+\ep}}+\frac{|\phi|^{p+1}+|\nabb\phi|^2}{r}\\
&\les E[\phi](\pa \cD_{u_1, u_2}^{-u_1})+\int_{\pa\cD_{u_1, u_2}^{-u_1}} r_+^{-2}|\phi|^2,\\
&\les  E[\phi]( \{t=0, 2|u_2|\leq |x|\})+E[\phi]( \H_{u_2}^{-u_1})+E[\phi](\Hb_{-u_1}^{u_1, u_2})\\
&\les (u_2)_+^{-\ga}\mathcal{E}_{\ga}[\phi]
\end{align*}
for all $u_1<u_2 <-1$. The last step follows from the above energy conservation in the exterior region and the definition of weighted energy $\mathcal{E}_{\ga}[\phi]$ for $\ga\geq 0$. By setting $u_1\rightarrow -\infty$, we obtain the integrated local energy decay estimates in the exterior region, which in particular implies estimate \eqref{eq:ILE:cD:thm} of the main theorem for $u<-1$.

Similarly in the interior region, from \eqref{eq:ILE:0:cD} and the above energy conservation, by using Hardy's inequality to bound the lower order term $r_+^{-2}|\phi|^2$, we obtain that
\begin{align*}
 & \iint_{\cD_{u_1, u_2}^{v_0}} \frac{|\pa\phi|^2+|r_+^{-1}\phi|^2}{(1+r)^{1+\ep}}+\frac{|\phi|^{p+1}+|\nabb\phi|^2}{r}\\
 &\les E[\phi](\pa \cD_{u_1, u_2}^{v_0})+\int_{\pa\cD_{u_1, u_2}^{v_0}} r_+^{-2}|\phi|^2,\\
&\les  E[\phi]( \Si_{u_1}^{v_0})+E[\phi]( \Si_{u_2}^{v_)})+E[\phi](\Hb_{v_0}^{u_1, u_2})+\int_{\Si_{u_1}^{v_0}}r_+^{-2}|\phi|^2\\
&\les E[\phi]( \Si_{u_1}^{v_0})+\int_{\Si_{u_1}^{v_0}}r_+^{-2}|\phi|^2
\end{align*}
for $-1\leq u_1<u_2$, $v_0\geq 2+2u_2$. Now letting $v_0\rightarrow \infty$ and using Hardy's inequality again, we obtain the integrated local energy decay estimates adapted to the foliation $\Si_u$ in the interior region
\begin{align}
  \label{eq:ILE:cD}
   E[\phi](\Si_{u_2})+\iint_{\cD_{u_1, u_2}} \frac{|\pa\phi|^2+|r_+^{-1}\phi|^2}{(1+r)^{1+\ep}}+\frac{|\phi|^{p+1}+|\nabb\phi|^2}{r}&\les E[\phi](\Si_{u_1} ).
\end{align}
Here we used the fact that the solution $\phi\rightarrow 0$ at the null infinity with finite energy initial data. Unlike the situation in the exterior region, the quantitative decay estimates will follow by combining these estimates with a type of $r$-weighted energy estimates, which will be introduced in the next section.

Finally in view of \eqref{eq:ILE:0:Rd}, by using the energy conservation as well as Hardy's inequality, we derive the integrated local energy decay estimate on the whole spacetime
\begin{align*}
   \int_0^t\int_{\mathbb{R}^d}\frac{|\pa\phi|^2+|r_+^{-1}\phi|^2}{(1+r)^{1+\ep}}+\frac{|\phi|^{p+1}+|\nabb\phi|^2}{r} &\les E[\phi]({0}\times\mathbb{R}^{d}). \end{align*}
Since the implicit constant is independent of $t$, the integrated local energy decay estimate \eqref{eq:ILE:1p5} follows by letting $t\rightarrow \infty$ in the second inequality.

\section{The $r$-weighted energy estimates and the quantitative flux decay}
The integrated local energy decay estimates obtained in the previous section are not sufficient to conclude the inverse polynomial decay of the energy flux \eqref{eq:ILE:cD:thm}. It has to be combined with a type of $r$-weighted energy estimates original introduced by Dafermos-Rodnianski in \cite{newapp}. In this section, according to the main theorem, we assume that
\[
\frac{d+1}{d-1}<p\leq \frac{d+2}{d-2},\quad 1<\ga_0<\min\{2, \f12 (d-1)(p-1)\}.
\]
In the energy identity \eqref{eq:energy:id}, choose the vector fields $X$, $Y$ and the function $\chi$ as follows:
\[
X=r^{\gamma} L,\quad Y=\frac{d-1}{4}\ga r^{\ga-2}|\phi|^2L, \quad \chi=\frac{d-1}{2}r^{\ga-1}
\]
for $1\leq \ga\leq \ga_0$.
We then compute that
\[
\nabla_{L}X=\gamma r^{\gamma-1}L,\quad \nabla_{\Lb}X= -\ga r^{\ga-1}L,\quad \nabla_{e_A}X=r^{\gamma-1} e_A.
\]
In particular the non-vanishing components of the deformation tensor $\pi_{\mu\nu}^X$ are
\[
\pi^X_{L\Lb}=-\gamma r^{\gamma-1},\quad \pi^X_{\Lb\Lb}=2\gamma r^{\gamma-1},\quad \pi^X_{e_A e_A}=r^{\ga-1}.
\]
Denote $\psi=r^{\frac{d-1}{2}}\phi$ and $c_d=\frac{(d-1)(d-3)}{4}$.
Then we can compute that
\begin{align*}
&div(Y)+T[\phi]^{\mu\nu}\pi^X_{\mu\nu}+
\chi \pa_\mu\phi \pa^\mu\phi  +\chi\phi\Box\phi-\f12 \Box\chi |\phi|^2\\
&=-\f12\ga r^{\gamma-1}(|\nabb\phi|^2+\frac{2}{p+1}|\phi|^{p+1})+r^{\ga-1}(|\nabb\phi|^2-\frac{d-1}{2}\pa^\mu \phi \pa_\mu\phi-\frac{d-1}{p+1}|\phi|^{p+1})\\
&\quad +\f12\ga r^{\ga-1}|L\phi|^2+
\chi \pa_\mu\phi \pa^\mu\phi  +\chi|\phi|^{p+1}-\f12 \Box\chi |\phi|^2+div(Y)\\
&=\f12 r^{\ga-1}(\ga|L\phi|^2+(2-\ga)|\nabb\phi|^2 )+(\frac{d-1}{2} -\frac{\ga  + d-1 }{p+1}) r^{\ga-1}|\phi|^{p+1}\\
&\quad -\frac{(d-1)(\ga-1)(d+\ga-3)}{4}r^{\ga-3}|\phi|^2+\frac{(d-1)\ga}{4}(L(r^{\ga-2}|\phi|^2)+(d-1)r^{\ga-3}|\phi|^2)\\
&=\f12 r^{\ga-d}(\ga|L\psi|^2+(2-\ga)(|\nabb\psi|^2+c_d r^{-2}|\psi|^2) )+(\frac{d-1}{2} -\frac{\ga  + d-1 }{p+1}) r^{\ga-1}|\phi|^{p+1} .
\end{align*}
Since $d\geq 3$ and $1\leq  \ga \leq \ga_0$, the quadratic part is nonnegative. On the other hand, as we have assumed $$\ga\leq \ga_0<\f12 (d-1)(p-1)$$
in this section,
we conclude that the coefficient of the potential part
\[
\frac{d-1}{2} -\frac{\ga  + d-1 }{p+1}=\frac{(d-1)(p-1)-2\ga}{2(p+1)}\geq \frac{\f12(d-1)(p-1)-\ga_0}{p+1}> 0
\]
also has a positive sign. This in particular means that for the above chosen vector fields $X$, $Y$ and function $\chi$, the bulk integral on the right hand side of the energy identity \eqref{eq:energy:id} is nonnegative.

Next we need to compute the left hand side of the energy identity \eqref{eq:energy:id}. By using Stokes' formula, we derive that
\begin{equation*}
\begin{split}
\iint_{\mathcal{D}}\pa^\mu J^{X, Y, \chi}_\mu[\phi]d\vol=\int_{\pa\mathcal{D}}i_{ J^{X, Y, \chi}[\phi]}d\vol.
\end{split}
\end{equation*}
 Here $i_Z \eta$ is the contraction of the differential form $\eta$ with the vector field $Z$.In this paper, the region $\mathcal{D}$ will always be the regular region bounded by the outgoing null hypersurface $\H_u$, the incoming null hypersurface $\Hb_v$ and the constant $t$-slice of the Minkowski space. In particular, the boundary $\pa\mathcal{D}$ consists these three kinds of hypersurfaces.
Recall the volume form
\[
d\vol=dtdx=2 r^{d-1} du d v d\om .
\]
Hence for the outgoing null hypersurface $\H_u$, we have
\begin{align*}
i_{J^{X, Y, \chi}[\phi]}d\vol&=2(J^{X, Y,\chi}[\phi])^{\Lb}r^{d-1}dvd\om= -( T[\phi]_{L \nu}X^\nu -
\f12 L\chi |\phi|^2 + \f12 \chi\cdot L |\phi|^2 +Y_L)  r^{d-1}dv d\om\\
&=-( |L\phi|^2 r^\ga -
\f12 L\chi |\phi|^2 +   \chi \phi\cdot L \phi )  r^{d-1}dv d\om.
\end{align*}
Recall that $\psi=r^{\frac{d-1}{2}}\phi$. We can write that
\begin{align*}
  |L\psi|^2=(r^{\frac{d-1}{2}}L\phi+\frac{d-1}{2}r^{\frac{d-3}{2}}\phi)^2=r^{d-1}|L\phi|^2+(d-1)r^{d-2}\phi L\phi+\frac{(d-1)^2}{4}r^{d-3}|\phi|^2.
\end{align*}
Now for the lower order terms, we show that
\begin{align*}
  &(d-1)r^{d-2} r^\ga \phi L\phi+r^\ga \frac{(d-1)^2}{4}r^{d-3}|\phi|^2+\f12 L\chi |\phi|^2 r^{d-1}-\chi r^{d-1}\phi L\phi \\
  &=((d-1)r^{d-2} r^\ga-\chi r^{d-1} )\phi L\phi+ \f12 (r^\ga \frac{(d-1)^2}{2}r^{-2} + L\chi   )r^{d-1}|\phi|^2\\
 &=\frac{d-1}{4}r^{d+\ga-2} L|\phi|^2+ \frac{d-1}{4} (d-1+\ga-1)  r^{d+\ga-3}|\phi|^2\\
 &=\frac{d-1}{4}L ( r^{d+\ga-2} |\phi|^2) .
\end{align*}
Thus on the outgoing null hypersurface, we have
\begin{align*}
  i_{J^{X, Y,\chi}[\phi]}d\vol
&=-  |L\psi|^2 r^\ga   dv d\om  + \frac{d-1}{4}L  (  r^{d+\ga-2} |\phi|^2)dvd\om.
\end{align*}
 The last term in the above identity will be cancelled by using integration by parts. We thus obtain a nonnegative weighted energy flux through the outgoing null hypersurface $\H_u$.

Next on the incoming null hypersurface $\Hb_v$, similarly we can show that
\begin{align*}
  i_{J^{X, Y, \chi}[\phi]}d\vol&=-2(J^{X, Y,\chi}[\phi])^{L}r^{d-1}dud\om= ( T[\phi]_{\Lb \nu}X^\nu -
\f12 \Lb\chi |\phi|^2 + \f12 \chi\cdot \Lb |\phi|^2 +Y_{\Lb})  r^{d-1}du d\om\\
&=(  (|\nabb\phi|^2+\frac{2}{p+1}|\phi|^{p+1})r^{\ga} -
\f12 \Lb\chi |\phi|^2 +   \chi \phi\cdot \Lb \phi -\f12 (d-1)\ga r^{\ga-2}|\phi|^2)  r^{d-1}du d\om.
\end{align*}
For the lower order terms, we show that
\begin{align*}
  & -\f12 \Lb\chi |\phi|^2 r^{d-1}+\chi r^{d-1}\phi \Lb\phi-\f12 (d-1)\ga r^{d+\ga-3}|\phi|^2 \\
  &=\frac{d-1}{4} r^{d+\ga-2}  \Lb|\phi|^2- \frac{(d-1)(\ga+1)}{4}r^{d+\ga-3}|\phi|^2\\
 &=\frac{d-1}{4}\Lb( r^{d+\ga-2} |\phi|^2)- \frac{d-1}{4}\Lb( r^{d+\ga-2} )|\phi|^2-\frac{(d-1)(\ga+1)}{4}r^{d+\ga-3}|\phi|^2 \\
 &=\frac{d-1}{4}\Lb (r^{d+\ga-2} |\phi|^2 )+\frac{(d-1)(d-3)}{4}r^{d+\ga-3}|\phi|^2.
\end{align*}
Thus on the incoming null hypersurface $\Hb_v$, we have
\begin{align*}
  i_{J^{X ,Y,\chi}[\phi]}d\vol
&= r^\ga (|\nabb\psi|^2+\frac{2}{p+1}|\phi|^{p+1}r^{d-1}+c_d r^{-2}|\psi|^2  )  du d\om  + \frac{d-1}{4}\Lb  (  r^{d+\ga-2} |\phi|^2)dud\om.
\end{align*}
Finally on the constant $t$-slice, we first have
\begin{align*}
  i_{J^{X, Y,\chi}[\phi]}d\vol&=(J^{X, Y,\chi}[\phi])^{0}dx= -( T[\phi]_{0 \nu}X^\nu -
\f12 \pa_t\chi |\phi|^2 + \f12 \chi\cdot \pa_t |\phi|^2 +Y_0)  dx\\
&=-\f12( r^\ga (|L\phi|^2+|\nabb\phi|^2+\frac{2}{p+1}|\phi|^{p+1})  +   \chi   \pa_t |\phi|^2-\frac{d-1}{2}\ga r^{\ga-2}|\phi|^2 )  dx.
\end{align*}
By writing the above integral in terms of the weighted function  $\psi=r^{\frac{d-1}{2}}\phi$, similar to the case on the outgoing null hypersurface, we can compute the lower order terms
\begin{align*}
  &(d-1)r^{d+\ga-2} \phi L\phi +  \frac{(d-1)^2}{4}r^{d+\ga-3}|\phi|^2 -\chi r^{d-1} \pa_t|\phi|^2+\frac{d-1}{2}\ga r^{d+\ga-3}|\phi|^2 \\
  &=(d-1)r^{d+\ga-2} \phi \pa_r\phi+ \frac{d-1}{4}(d+2\ga-1 )r^{d+\ga-3}|\phi|^2\\
 &=\frac{d-1}{2}\pa_r ( r^{d+\ga-2} |\phi|^2)-\frac{(d-1)(d-3)}{4}r^{d+\ga-3} |\phi|^2.
\end{align*}
Thus on the constant $t$-slice $\Si_t$, we have
\begin{align*}
  i_{J^{X,Y,\chi}[\phi]}d\vol
&=-\f12 r^\ga(  |L\psi|^2 +|\nabb\psi|^2+\frac{2}{p+1}|\phi|^{p+1}r^{d-1}+c_d r^{-2}|\psi|^2) dr d\om\\
& \quad + \frac{d-1}{4}\pa_r  (  r^{d+\ga-2} |\phi|^2 )drd\om.
\end{align*}
Based on these computations, we are now ready to derive the necessary $r$-weighted energy estimates for the solution $\phi$.
In the exterior region, take the region $\cD$ to be $\cD_{u_1, u_2}^{-u_1}$ with $u_1\leq u_2\leq -1$. Note that this region is bounded by the initial hypersurface, the outgoing null hypersurface $\H_{u_2}^{-u_1}$ and the incoming null hypersurface $\Hb_{-u_1}^{u_1, u_2}$. On the boundary $\pa\cD_{u_1, u_2}^{-u_1}$, since
\begin{align*}
  \int_{\{t=0, |u_2|\leq \frac{|x|}{2}\leq |u_1|\}} \pa_r  (  r^{d+\ga-2} |\phi|^2 )drd\om+\int_{\Hb_{-u_1}^{u_1, u_2}} \Lb(r^{d+\ga-2} |\phi|^2)dud\om- \int_{\H_{u_2}^{-u_1}}L(r^{d+\ga-2} |\phi|^2)dvd\om=0,
\end{align*}
the above computations together with the energy identity \eqref{eq:energy:id} then lead to the following weighted energy identity in the exterior region
\begin{align*}
  &\iint_{\cD_{u_1, u_2}^{-u_1}}  r^{\ga-d}(\ga|L\psi|^2+(2-\ga)(|\nabb\psi|^2+c_d r^{-2}|\psi|^2) )+2(\frac{d-1}{2} -\frac{\ga  + d-1 }{p+1}) r^{\ga-1}|\phi|^{p+1}\\
  &+\int_{\H_{u_2}^{-u_1}}2r^{\ga}|L\psi|^2dvd\om+\int_{\Hb_{-u_1}^{u_1, u_2}}2 r^\ga (|\nabb\psi|^2+\frac{2}{p+1}|\phi|^{p+1}r^{d-1}+c_d r^{-2}|\psi|^2  )  du d\om\\
  &=\int_{\{t=0, |u_2|\leq \frac{|x|}{2}\leq |u_1|\}}  r^\ga(  |L\psi|^2 +|\nabb\psi|^2+\frac{2}{p+1}|\phi|^{p+1}r^{d-1}+c_d r^{-2}|\psi|^2) dr d\om
\end{align*}
for all $1\leq \ga \leq \ga_0$ and $u_1<u_2\leq -1$. Here recall that $c_d=\frac{(d-1)(d-3)}{4}$ and $\psi=r^{\frac{d-1}{2}}\phi$.

By setting $\ga=\ga_0$, the right hand side can be bounded by the initial weighted energy flux $\mathcal{E}_{\ga_0}[\phi]$ together with Hardy's inequality for controlling the integral of $r^{-2}|\phi|^2$. Since the left hand side is nonnegative due to the assumption on $p$ and $\ga_0$, the above identity in particular leads to the following weighted energy estimate in the exterior region
\begin{equation}
\label{eq:PWE:ex}
  \begin{split}
    \int_{-\infty}^{u_2}\int_{\H_u} r^{\ga_0-d}( |L\psi|^2+ |\nabb\psi|^2)+  r^{\ga_0-1}|\phi|^{p+1} +\int_{\H_{u_2} } r^{\ga_0}|L\psi|^2dvd\om \les \mathcal{E}_{\ga_0}[\phi]
  \end{split}
\end{equation}
for all $u_2\leq -1$ by letting $u_1\rightarrow -\infty$. Here the implicit constant relies only on $\ga_0$, $p$ and $d$. In particular, the $r$-weighted energy estimate \eqref{eq:PEW:thm} of the main theorem holds when the integral is restricted to the exterior region by noting that $v_+\les r$ in this region.

For the case in the interior region, consider the domain $\cD_{u_1, u_2}^{v_0}$ with $-1<u_1<u_2$ and $v_0\geq 2u_2+2$. Similarly, we can  obtain the following weighted energy identity
\begin{align*}
  &\iint_{\cD_{u_1, u_2}^{v_0}}  r^{\ga-d}(\ga|L\psi|^2+(2-\ga)(|\nabb\psi|^2+c_d r^{-2}|\psi|^2) )+2(\frac{d-1}{2} -\frac{\ga  + d-1 }{p+1}) r^{\ga-1}|\phi|^{p+1}\\
  &+\int_{\{t=2+2u_2, |x|\leq 2\}}  r^\ga(  |L\psi|^2 +|\nabb\psi|^2+\frac{2}{p+1}|\phi|^{p+1}r^{d-1}+c_d r^{-2}|\psi|^2) dr d\om\\
  &+\int_{\H_{u_2}^{v_0}}2r^{\ga}|L\psi|^2dvd\om+\int_{\Hb_{v_0}^{u_1, u_2}}2 r^\ga (|\nabb\psi|^2+\frac{2}{p+1}|\phi|^{p+1}r^{d-1}+c_d r^{-2}|\psi|^2  )  du d\om\\
  &= \int_{\{t=2+2u_1, |x|\leq 2\}}  r^\ga(  |L\psi|^2 +|\nabb\psi|^2+\frac{2}{p+1}|\phi|^{p+1}r^{d-1}+c_d r^{-2}|\psi|^2) dr d\om +\int_{\H_{u_1}^{v_0}}2r^{\ga}|L\psi|^2dvd\om.
\end{align*}
Again the left hand side is nonnegative. The first integral on the right hand side can be trivially bounded by the energy flux as $|x|\leq 2$ up to a constant. Thus by letting $v_0\rightarrow \infty$, we conclude from the previous identity that
\begin{align}
\notag
  &\int_{u_1}^{u_2}\int_{\Si_{u}}  r^{\ga-1}(r^{1-d}|L\psi|^2+|\nabb\phi|^2+|\phi|^{p+1}+c_d r^{-2}|\phi|^2) +\int_{\H_{u_2} }r^{\ga}|L\psi|^2dvd\om \\
  \label{eq:PEW:ga:general}
  & \les E[\phi](\Si_{u_1}) +\int_{\H_{u_1} }r^{\ga}|L\psi|^2dvd\om
\end{align}
for all $1\leq \ga \leq \ga_0$. Setting $u_1=-1$ and $\ga=\ga_0$, in view of the weighted energy decay estimate \eqref{eq:PWE:ex} and the energy bound \eqref{eq:Eflux:decay:ex} in the exterior region, we in particular derive from the above inequality that
\begin{align*}
   &\int_{-1}^{u_2}\int_{\H_{u}}  r^{\ga_0-d}|L\psi|^2 +\int_{\H_{u_2} }r^{\ga_0}|L\psi|^2dvd\om \les E[\phi](\Si_{-1}) +\int_{\H_{-1} }r^{\ga_0}|L\psi|^2dvd\om \les \mathcal{E}_{\ga_0}[\phi]
\end{align*}
for all $-1\leq u_2$, which in particular implies that we can extract a dyadic sequence $\{u_k\}_{3}^{\infty}$ such that
\begin{align*}
  u_3=1,\quad  2 u_k&\leq u_{k+1}\leq  \La u_k,\\
   \int_{\H_{u_k}}r^{\ga_0-1}|L\psi|^2 dvd\om &\leq C (1+u_k)^{-1}\mathcal{E}_{\ga_0}[\phi]
\end{align*}
for some constants $\La$ and $C$ depending only on $p$, $d$ and $\ga_0$. Since the previous inequality in particular implies that
$$\int_{\H_u}r^{\ga_0}|L\psi|^2dvd\om\les \mathcal{E}_{\ga_0}[\phi],\quad \forall -1\leq u,$$
interpolation then leads to
\begin{equation}
\label{eq:PEW:r1}
  \int_{\H_{u_k}}r |L\psi|^2 dvd\om \les u_k^{1-\ga_0}\mathcal{E}_{\ga_0}[\phi],\quad k\geq 3.
\end{equation}
Here recall that we have assumed that $\ga_0>1$.

Now we note that
\begin{align*}
    \int_{u_1}^{u_2} \int_{\H_u} |L\phi|^2
  &=\int_{u_1}^{u_2}\int_{\H_u} r^{1-d}|L\psi|^2 -(d-1)r^{-1}\phi L\phi-\frac{(d-1)^2}{4}r^{-2}|\phi|^2\\
  &=\int_{u_1}^{u_2}\int_{\H_u} r^{1-d}|L\psi|^2-(d-1) \int_{u_1}^{u_2}\int_{\H_u} L(r^{d-2}|\phi|^2)+\frac{3-d}{2}r^{d-3}|\phi|^2 dvd\om du\\
  &\leq \int_{u_1}^{u_2}\int_{\H_u} r^{1-d}|L\psi|^2+c_d r^{-2}|\phi|^2+(d-1)\int_{2u_1+2}^{2u_2+2}\int_{|x|=2}r^{d-2}|\phi|^2 d\om dt\\
  &\les \int_{u_1}^{u_2}\int_{\H_u} r^{1-d}|L\psi|^2+c_d r^{-2}|\phi|^2+\int_{2u_1+2}^{2u_2+2}\int_{|x|\leq 2} |\pa\phi|^2+r_+^{-2}|\phi|^2dxdt
\end{align*}
Thus by letting  $\ga=1$ in \eqref{eq:PEW:ga:general} and  in view of the integrated local energy decay estimate \eqref{eq:ILE:cD}, we can show that
\begin{align*}
  \int_{u_1}^{u_2} E[\phi](\Si_u)du&=\int_{u_1}^{u_2}(\int_{|x|\leq 2}|\pa\phi|^2+\frac{2}{p+1}|\phi|^{p+1} +\int_{\H_u} |L\phi|^2+|\nabb\phi|^2+\frac{2}{p+1}|\phi|^{p+1} )du\\
  &\les \int_{u_1}^{u_2}\int_{\Si_u}\frac{|\pa\phi|^2 +r_+^{-2}|\phi|^2}{(1+r)^{1+\ep}}+  r^{1-d}|L\psi|^2+|\nabb\phi|^2+\frac{2}{p+1}|\phi|^{p+1}+c_d r^{-2}|\phi|^2\\
  &\les E[\phi](\Si_{u_1})+\int_{\H_{u_1}}r|L\psi|^2dvd\om.
\end{align*}
As the energy flux $E[\phi](\Si_u)$ is non-increasing in view of \eqref{eq:ILE:cD}, by setting $u_1=-1$, we conclude that
\begin{align*}
   \f12 u E[\phi](\Si_{2u})\les \int_{\f12 u}^{u}E[\phi](\Si_{u'})du'\les E[\phi](\Si_{-1})+\int_{\H_{-1}}r|L\psi|^2 dvd\om\les \mathcal{E}_{\ga_0}[\phi]
\end{align*}
for all $1<u$. Thus we derive a weak decay estimate for the energy flux
\begin{equation*}
  E[\phi](\Si_{u})\les u_+^{-1} \mathcal{E}_{\ga_0}[\phi],\quad \forall -1\leq u.
\end{equation*}
By using this weak decay estimate and in view of \eqref{eq:PEW:r1}, we derive from the previous estimate that
\begin{align*}
  (u_{k+1}-u_k)E[\phi](\Si_{u_{k+1}})\leq \int_{u_k}^{u_{k+1}}E[\phi](\Si_u)du\les E[\phi](\Si_{u_k})+\int_{\H_{u_k}}r|L\psi|^2dvd\om\les u_k^{1-\ga_0}\mathcal{E}_{\ga_0}[\phi].
\end{align*}
Here we used the assumption $1<\ga_0\leq \ga_0<2$. Since the sequence $u_{k}$ is dyadic and verifies the relation $2u_k\leq u_{k+1}\leq \La u_{k}$, we therefore can demonstrate that
\begin{align*}
 E[\phi](u) \leq E[\phi](u_{k+1})\les (u_{k+1}-u_k)^{-1} u_k^{1-\ga_0}\mathcal{E}_{\ga_0}[\phi]\les u_{k+1}^{-\ga_0}\mathcal{E}_{\ga_0}[\phi]\les u^{-\ga_0} \mathcal{E}_{\ga_0}[\phi]
\end{align*}
for all $u\in [u_{k+1}, u_{k+2}]$, $k\geq 3$. This energy decay estimate together with the integrated local energy decay estimate \eqref{eq:ILE:cD} leads to the quantitative inverse polynomial decay estimate \eqref{eq:ILE:cD:thm} of the main theorem.

As for the $r$-weighted energy estimate \eqref{eq:PEW:thm} of the main theorem in the interior region, from \eqref{eq:PEW:ga:general}, we in particular derive that
\begin{align*}
  &\int_{-1}^{\infty}\int_{\Si_{u}}  r^{\ga_0-d}|L\psi|^2+r^{\ga_0-1}|\nabb\phi|^2+r^{\ga_0-1}|\phi|^{p+1}
  \les E[\phi](\Si_{-1}) +\int_{\H_{-1} }r^{\ga_0}|L\psi|^2dvd\om\les \mathcal{E}_{\ga_0}[\phi].
\end{align*}
Thus to conclude \eqref{eq:PEW:thm}, it remains to improve the bound for the potential part. Note that by using the energy flux decay estimate \eqref{eq:ILE:cD:thm}, we can show that
\begin{align*}
  \int_{-1}^{\infty}\int_{\Si_{u}}  v_+^{\ga_0-1-\ep}|\phi|^{p+1}
  &\leq  \int_{-1}^{\infty}\int_{\Si_{u}}  (r^{\ga_0-1-\ep}+u_+^{\ga_0-1-\ep})|\phi|^{p+1}\\
  &\les  \mathcal{E}_{\ga_0}[\phi]+\int_{-1}^{\infty} u_+^{\ga_0-1-\ep} E[\phi](\Si_u)du\\
  &\les \mathcal{E}_{\ga_0}[\phi]+\int_{-1}^{\infty} u_+^{\ga_0-1-\ep} u_+^{-\ga_0}\mathcal{E}_{\ga_0}[\phi]du\\
  &\les \mathcal{E}_{\ga_0}[\phi].
\end{align*}
Hence we finished the proof for the $r$-weighted energy estimates \eqref{eq:PEW:thm} of the main theorem.

\section{Proof for the scattering results}
In this section, we prove our scattering results. Assume that $p$, $d$ and $\ga_0$ verifies the following relation
$$\max\{\frac{4}{p-1}-d+2, 1\}<\ga_0<\min\{\f12 (p-1)(d-1), 2\},$$
which implies that
\[
p>p(d)=\frac{1+\sqrt{d^2+4d-4}}{d-1}.
\]
It is clear that $p(d)<p^*(d)=\frac{d+2+\sqrt{d^2+8d}}{2(d-1)}<\frac{d+3}{d-1}$, where $p^*(d)$ is the lower bound in the work of Ginibre-Velo for dimension $d=3$, $4$.

We first prove that our scattering result follows from the uniform spacetime bound \eqref{eq:bd4:Lq:spacebd}.
\begin{Prop}
  \label{prop:scattering:finiteLpnorm}
  Let $\phi$ be solution of \eqref{eq:NLW:semi:alld} with initial data $(\phi_0, \phi_1)\in (\dot{H}_x^1\cap \dot H^{s_p}_x)\times (L^2_x\cap \dot{H}^{s_p-1}_x)$. Suppose that the spacetime norm $\|\phi\|_{L_{t, x}^{\frac{(d+1)(p-1)}{2}}}$ is bounded. Then the solution $\phi$ scatters in $\dot{H}^{s}\times \dot{H}^{s-1}$ for all $s_p\leq s\leq 1$, that is , there exists $(\phi_0^{\pm}, \phi_1^{\pm})$ such that
  \begin{align*}
    \lim\limits_{t\rightarrow\pm\infty}\|(\phi(t, x),\pa_t\phi(t,x))-\mathbf{L}(t)(\phi_0^{\pm}(x), \phi_1^{\pm}(x))\|_{\dot{H}_x^s\times \dot{H}_x^{s-1}}=0.
  \end{align*}
\end{Prop}
The proof is inspired by that in \cite{Shen17:3DNLW:scatter} for the case when $d=3$. For readers' interests, we repeat the proof here.
\begin{proof}
  Since the wave equation is time invertible, it suffices to prove the scattering result in the future direction. Moreover, by interpolation, we only need to prove the solution scatters in the endpoint case when $s=s_p$ and $s=1$.

  For the case when $s=s_p$ and $1<p\leq 1+\frac{4}{d-1}$, using Strichartz estimate, we show that
  \begin{align*}
    &\|\mathbf{L}(-t_2)(\phi(t_2, x), \pa_t\phi(t_2, x))-\mathbf{L}(-t_1)(\phi(t_1, x), \pa_t\phi(t_1, x))\|_{\dot{H}_x^{s_p}\times\dot{H}_x^{s_p-1}}\\
    &=\|(\phi(t_2, x), \pa_t\phi(t_2, x))-\mathbf{L}(t_2-t_1)(\phi(t_1, x), \pa_t\phi(t_1, x))\|_{\dot{H}_x^{s_p}\times\dot{H}_x^{s_p-1}}\\
    &\leq C_{d, p} \||\phi|^{p-1}\phi\|_{L_{t,x}^{\frac{(d+1)(p-1)}{2p}}([t_1, t_2]\times\mathbb{R}^d)} = C_{d, p} \|\phi\|_{L_{t,x}^{\frac{(d+1)(p-1)}{2}}([t_1, t_2]\times\mathbb{R}^d)}^p.
  \end{align*}
  Here the constant $C_{d, p}$ replies only on $d$ and $p$. The restriction $1<p\leq 1+\frac{4}{d-1}$ on $p$ is to guarantee the pairs used in the above Strichartz estimate are admissible. As $\|\phi\|_{L_{t, x}^{\frac{(d+1)(p-1)}{2}}}$ is finite, we conclude that $\mathbf{L}(-t)(\phi(t, x), \pa_t\phi(t, x))$ converges to some pair $(\phi_0^+, \phi_1^+)$ in $\dot{H}^{s_p}\times\dot{H}^{s_p-1}$ as $t\rightarrow +\infty$. In particular,
  \begin{align*}
    \lim\limits_{t\rightarrow +\infty}\|(\phi(t, x),\pa_t\phi(t,x))-\mathbf{L}(t)(\phi_0^{+}(x), \phi_1^{+}(x))\|_{\dot{H}_x^{s_p}\times \dot{H}_x^{s_p-1}}=0.
  \end{align*}
  For the case when $s=s_p$ and $\frac{4}{d-1}<p-1\leq \frac{4}{d-2}$ or $s=1$, note that we in particular have $s\geq \frac{1}{2}$. On any finite time interval $[t_1, t_2]$, applying the Strichartz estimate to the fractional derivatives $\nabla_x^{s-\frac{1}{2}}\phi$ of $\phi$, we can estimate that
  \begin{equation}
  \label{eq:bd4:nablaf12phi}
  \begin{split}
    &\|\nabla_x^{s-\frac{1}{2}}\phi\|_{L_{t, x}^{\frac{2(d+1)}{d-1}}([t_1, t_2]\times\mathbb{R}^d)}\\
    &\leq C_{d}(\|\phi(t_1, x)\|_{\dot{H}_x^{s}}+\|\pa_t\phi(t_1, x)\|_{\dot{H}_x^{s-1}}+\|\nabla_x^{s-\f12}(|\phi|^{p-1}\phi)\|_{L_{t, x}^{\frac{2(d+1)}{d+3}}([t_1, t_2]\times\mathbb{R}^d)})\\
    &\leq C_{d}(\|\phi(t_1, x)\|_{\dot{H}_x^{s}}+\|\pa_t\phi(t_1, x)\|_{\dot{H}_x^{s-1}}+\|\nabla_x^{s-\f12} \phi\|_{L_{t, x}^{\frac{2(d+1)}{d-1}}([t_1, t_2]\times\mathbb{R}^d)}\|\phi\|^{p-1}_{L_{t, x}^{\frac{(d+1)(p-1)}{2}}([t_1, t_2]\times\mathbb{R}^d)})
    \end{split}
  \end{equation}
  for some constant $C_{d}$ relying only on $d$. As $\|\phi\|_{L_{t, x}^{\frac{(d+1)(p-1)}{2}}}$ is bounded, take $t_1$ large enough such that
  \[
  C_{d}\|\phi\|^{p-1}_{L_{t, x}^{\frac{(d+1)(p-1)}{2}}([t_1, \infty)\times\mathbb{R}^d)}<\frac{1}{2}.
  \]
  We conclude from the previous inequality that
  \begin{align*}
    \|\nabla_x^{s-\frac{1}{2}}\phi\|_{L_{t, x}^{\frac{2(d+1)}{d-1}}([t_1, \infty)\times\mathbb{R}^d)}<\infty.
  \end{align*}
  Here the boundedness of $\|(\phi(t_1, x), \pa_t\phi(t_1, x))\|_{\dot{H}_x^1\times L_x^2}$ follows from the energy conservation and then together with Strichartz estimate leads to the finiteness of $\|(\phi(t_1, x), \pa_t\phi(t_1, x))\|_{\dot{H}_x^{s_p}\times \dot{H}_x^{s_p-1}}$.

  Now we run the same argument to estimate that
    \begin{align*}
    &\|(\phi(t_3, x), \pa_t\phi(t_3, x))-\mathbf{L}(t_3-t_2)(\phi(t_2, x), \pa_t\phi(t_2, x))\|_{\dot{H}_x^{s}\times\dot{H}_x^{s-1}}\\
    &\leq C_{d, p} \|\nabla_x^{s-\f12}(|\phi|^{p-1}\phi)\|_{L_{t,x}^{\frac{2(d+1)}{d+3}}([t_2, t_3]\times\mathbb{R}^d)} \\
    &\leq C_{d, p} \|\nabla_x^{s-\f12} \phi\|_{L_{t, x}^{\frac{2(d+1)}{d-1}}([t_2, t_3]\times\mathbb{R}^d)}\|\phi\|^{p-1}_{L_{t, x}^{\frac{(d+1)(p-1)}{2}}([t_2, t_3]\times\mathbb{R}^d)}
  \end{align*}
  for all $t_1<t_2<t_3$. This leads to the claim of the Proposition.
\end{proof}

Now to conclude our main theorem, it remains to prove the spacetime bound \eqref{eq:bd4:Lq:spacebd} restricted to the future as the bound in the past can be obtained in the same way. For the sub-conformal case $p(d)<p\leq \frac{d+3}{d-1}$, note that
\begin{align*}
  \frac{(d+1)(p-1)}{2}\leq p+1.
\end{align*}
By using the $r$-weighted energy estimate \eqref{eq:PEW:thm}, we conclude that
\begin{align*}
  \|\phi\|_{L_{t,x}^{\frac{(d+1)(p-1)}{2}}}\leq \|\phi v_+^{\frac{\ga_0-1-\ep}{p+1}}\|_{L_{t,x}^{p+1}}\| v_+^{-\frac{\ga_0-1-\ep}{p+1}}\|_{L_{t,x}^q}
\end{align*}
with
$$\frac{1}{q}+\frac{1}{p+1}=\frac{2}{(d+1)(p-1)}.$$
Since $\ga_0>\frac{4}{p-1}-d+2$, we in particular have
\begin{align*}
  \frac{\ga_0-1}{(p+1)(d+1)}>\frac{2}{(p-1)(d+1)}-\frac{1}{p+1}.
\end{align*}
Therefore
\begin{align*}
  \frac{\ga_0-1}{p+1}q>d+1.
\end{align*}
Choose $\ep$ sufficiently small such that
\begin{align*}
  \frac{\ga_0-1-\ep}{p+1}q>d+1,
\end{align*}
which implies that $\| v_+^{-\frac{\ga_0-1-\ep}{p+1}}\|_{L_{t,x}^q}$ is finite. We thus conclude that
\begin{align*}
  \|\phi\|_{L_{t,x}^{\frac{(d+1)(p-1)}{2}}}\les (\mathcal{E}_{\ga_0}[\phi])^{\frac{1}{p+1}}.
\end{align*}
This proves the bound \eqref{eq:bd4:Lq:spacebd} for the sub-conformal case.

Finally we prove the bound \eqref{eq:bd4:Lq:spacebd} for the sup-conformal case $\frac{d+3}{d-1}<p<\frac{d+2}{d-2}$, which implies that $1<\ga_0<2$. By choosing $0<\ep<\ga_0-1$, we derive from the $r$-weighted energy estimate \eqref{eq:PEW:thm} that
\begin{align*}
  \|\phi (1+t)^{\frac{\ga_0-1-\ep}{p+1}}\|_{L_{t, x}^{p+1}}\les (\mathcal{E}_{\ga_0}[\phi])^{\frac{1}{p+1}}.
\end{align*}
Here note that we are only interested in the estimates in the future and $v_+\geq \frac{1+t}{2}$.
On the other hand, from the energy conservation, we have
\begin{align*}
  \|\phi\|_{L_{t}^{\infty}L_{ x}^{p+1}}\les (\mathcal{E}_{\ga_0}[\phi])^{\frac{1}{p+1}}.
\end{align*}
Therefore for all $s\geq p+1$, we have the uniform mixed norm bound
\begin{align}
\label{eq:bd4LsLp1}
  \|\phi\|_{L_{t}^{s}([t_1, \infty))L_{ x}^{p+1}}\les (1+t_1)^{-\frac{\ga_0-1-\ep}{s}}(\mathcal{E}_{\ga_0}[\phi])^{\frac{1}{p+1}},\quad \forall s\geq p+1,\quad t_1\geq 0.
\end{align}
By using Sobolev embedding, we have
\begin{align*}
  \|\phi\|_{L_t^{\frac{2(d+1)}{d-1}}L_x^{q}}\les \|\nabla_x^{\f12}\phi\|_{L_{t,x}^{\frac{2(d+1)}{d-1}}},\quad \textnormal{ with } \frac{1}{q}+\frac{1}{2d}=\frac{d-1}{2(d+1)}.
\end{align*}
Since $\frac{d+3}{d-1}<p<\frac{d+2}{d-2}$, we can show that
\begin{align*}
  p+1<\frac{(d+1)(p-1)}{2}<q.
\end{align*}
Interpolation then implies that
\begin{align}
\label{eq:final}
  \|\phi\|_{L_{t,x}^{\frac{(d+1)(p-1)}{2}}}\leq \|\phi\|_{L_t^{q_1} L_{x}^{p+1}}^{\theta} \|\phi\|_{L_t^{\frac{2(d+1)}{d-1}}L_x^{q}}^{1-\theta}\les  \|\phi\|_{L_t^{q_1} L_{x}^{p+1}}^{\theta} \|\nabla_x^{\f12}\phi\|_{L_{t,x}^{\frac{2(d+1)}{d-1}}}^{1-\theta},
\end{align}
where
\begin{align*}
  \frac{\theta}{q_1}+\frac{(d-1)(1-\theta)}{2(d+1)}=\frac{2}{(d+1)(p-1)}=\frac{\theta}{p+1}+\frac{1-\theta}{q},\quad 0<\theta<1.
\end{align*}
Now the estimate \eqref{eq:bd4:nablaf12phi} in the proof of the above Proposition with $s=1$ then implies that
  \begin{equation*}
  \begin{split}
    \|\nabla_x^{\frac{1}{2}}\phi\|_{L_{t, x}^{\frac{2(d+1)}{d-1}}([t_1, t_2]\times\mathbb{R}^d)}
    &\leq C_{d}(\mathcal{E}_{0}[\phi]+  \|\phi\|_{L_t^{q_1}([t_1, t_2]) L_{x}^{p+1}}^{(p-1)\theta}\|\nabla_x^{\f12}\phi\|_{L_{t,x}^{\frac{2(d+1)}{d-1}}([t_1, t_2]\times\mathbb{R}^d)}^{1+(1-\theta)(p-1)}).
    \end{split}
  \end{equation*}
Note that $\frac{(d+1)(p-1)}{2}>\frac{2(d+1)}{d-1}$ for the sup-conformal case. We in particular conclude that
\[
p+1<\frac{(d+1)(p-1)}{2}<q_1<\infty.
\]
Thus from the decay estimate \eqref{eq:bd4LsLp1}, we can choose $t_1$ sufficiently large, depending on $\mathcal{E}_{\ga_0}[\phi]$, $p$, $\ga_0$ and  $d$, such that
\begin{align*}
  C_{d}\|\phi\|_{L_t^{q_1}([t_1, \infty)) L_{x}^{p+1}}^{(p-1)\theta}<(1+C_{d}\mathcal{E}_{\ga_0}[\phi])^{-(p-1)(1-\theta)},
\end{align*}
which guarantees that the function
\begin{align*}
  f(s)=C_{d}\mathcal{E}_{0}[\phi]+C_{d}\|\phi\|_{L_t^{q_1}([t_1, \infty)) L_{x}^{p+1}}^{(p-1)\theta} s^{1+(p-1)(1-\theta)}-s,\quad s\geq 0
\end{align*}
has exactly two distinct zeros. Therefore from the previous inequality and the continuity of the mixed norm $\|\nabla_x^{\frac{1}{2}}\phi\|_{L_{t, x}^{\frac{2(d+1)}{d-1}}([t_1, t_2]\times\mathbb{R}^d)}$, we derive that
\begin{align*}
  \|\nabla_x^{\frac{1}{2}}\phi\|_{L_{t, x}^{\frac{2(d+1)}{d-1}}([t_1, t_2]\times\mathbb{R}^d)}\les \mathcal{E}_{\ga_0}[\phi],\quad \forall t_2>t_1.
\end{align*}
In particular, the norm $\|\nabla_x^{\frac{1}{2}}\phi\|_{L_{t, x}^{\frac{2(d+1)}{d-1}} } $ is finite when restricted to the time interval $[t_1, \infty)$ for some large time $t_1>0$. For the norm on the finite interval $[0, t_1]$, we can divide this finite interval into small intervals $[t_{k}, t_{k+1}]$, on which, by using the energy conservation
\begin{align*}
  \|\phi\|_{L_{t}^{q_1}([t_k, t_{k+1}])L_x^{p+1}}\les (t_{k+1}-t_k)^{\frac{1}{q_1}}(\mathcal{E}_{\ga_0}[\phi])^{\frac{1}{p+1}}.
\end{align*}
Thus for the same reason, we can demonstrate that $\|\nabla_x^{\frac{1}{2}}\phi\|_{L_{t, x}^{\frac{2(d+1)}{d-1}} } $ is finite on these small intervals. This shows that the spacetime norm $\|\nabla_x^{\frac{1}{2}}\phi\|_{L_{t, x}^{\frac{2(d+1)}{d-1}} } $ is bounded, which leads to the claim \eqref{eq:bd4:Lq:spacebd} of the main theorem in view of the inequality \eqref{eq:final}.

\bibliography{shiwu}{}
\bibliographystyle{plain}


Beijing International Center for Mathematical Research, Peking University,
Beijing, China

\textsl{Email address}: shiwuyang@math.pku.edu.cn

\end{document}